\documentclass[12pt]{amsart}
\usepackage{amssymb}
\usepackage{amsxtra}
\usepackage{graphics}
\usepackage{charter}

\usepackage[margin=1in]{geometry}

\def\cI{\mathcal{I}}

\renewcommand{\setminus}{-}  
\renewcommand{\epsilon}{\varepsilon}  
\newcommand{\vecalpha}{\vec{\alpha}}

\newcommand{\norm}{{\mathbf N}}
\newcommand{\bad}{{\operatorname{bad}}}

\newcommand{\F}{{\mathbb F}}

\newcommand{\N}{{\mathbb N}}
\newcommand{\Q}{{\mathbb Q}}
\newcommand{\R}{{\mathbb R}}
\newcommand{\Z}{{\mathbb Z}}

\newcommand{\be}{\begin{enumerate}}
\newcommand{\ee}{\end{enumerate}}
\newcommand{\pp}{{\mathfrak p}}
\newcommand{\bb}{{\mathfrak b}}
\newcommand{\qq}{{\mathfrak q}}

\newcommand{\dd}{{\mathfrak d}}

\newcommand{\calA}{{\mathcal A}}

\newcommand{\calE}{{\mathcal E}}

\newcommand{\calL}{{\mathcal L}}
\newcommand{\calM}{{\mathcal M}}

\newcommand{\calO}{{\mathcal O}}
\newcommand{\calP}{{\mathcal P}}
\newcommand{\calQ}{{\mathcal Q}}

\newcommand{\calS}{{\mathcal S}}
\newcommand{\calT}{{\mathcal T}}

\newcommand{\calW}{{\mathcal W}}

\newcommand{\calZ}{{\mathcal Z}}

\newcommand{\OO}{{\mathcal O}}

\DeclareMathOperator{\ord}{ord}

\newcommand{\comment}[1]{}

\newcommand{\tors}{{\operatorname{tors}}}

\newcommand{\isom}{\simeq}

\newcommand{\Union}{\bigcup}
\newcommand{\union}{\cup}

\newtheorem{theorem}{Theorem}[section]
\newtheorem{lemma}[theorem]{Lemma}
\newtheorem{corollary}[theorem]{Corollary}
\newtheorem{proposition}[theorem]{Proposition}

\theoremstyle{definition}
\newtheorem{definition}[theorem]{Definition}
\newtheorem{question}[theorem]{Question}
\newtheorem{questions}[theorem]{Questions}

\theoremstyle{remark}
\newtheorem{remark}[theorem]{Remark}

\newtheorem{notation}[theorem]{Notation}

\begin{document}
\title[Hilbert's Tenth Problem and
Mazur's Conjectures]{Hilbert's Tenth Problem and Mazur's Conjectures in
  Complementary Subrings of Number Fields}
\author[K.\ Eisentr\"ager]{Kirsten Eisentr\"ager}
\address{Department of Mathematics, The Pennsylvania State University,
University Park, PA 16802, USA.}
\email{eisentra@math.psu.edu}
\urladdr{}
\author[G.\ Everest]{Graham Everest}
\address{School of Mathematics, University of East Anglia,
Norwich NR4 7TJ, UK.}
\email{g.everest@uea.ac.uk}
\author[A.\ Shlapentokh]{Alexandra Shlapentokh}
\address{Department of Mathematics,
East Carolina University,
Greenville, NC 27858, USA.}
\email{shlapentokha@ecu.edu}

\begin{abstract}
  We show that Hilbert's Tenth Problem is undecidable for
  complementary subrings of number fields and that the $\pp$-adic and
  archimedean ring versions of Mazur's conjectures do not hold in
  these rings.  More specifically, given a number field $K$, a
  positive integer $t>1$, and $t$ nonnegative computable real numbers
  $\delta_1,\ldots,\delta_t$ whose sum is one, we prove that the
  nonarchimedean primes of $K$ can be partitioned into $t$ disjoint
  recursive subsets $S_1, \dots, S_t$ of densities
  $\delta_1,\ldots,\delta_t$, respectively such that Hilbert's Tenth
  Problem is undecidable for each corresponding ring $\OO_{K,S_i}$. We
  also show that we can find a partition as above such that each ring
  $\OO_{K,S_i}$ possesses an infinite Diophantine set which is
  discrete in every topology of the field. The only assumption on $K$
  we need is that there is an elliptic curve of rank one defined over
  $K$.  \thanks{K.E.\ was partially supported by National Science
    Foundation grants DMS-0801123, DMS-1056703 and a Sloan Research
    Fellowship. A.S.\ was partially supported by National Science
    Foundation grant DMS-0650927 and a grant from the John Templeton
    Foundation.
    \\
    \indent Sadly, the second author passed away before the final
    version of the paper was completed.}
\end{abstract}

\maketitle

\section{Introduction}
Hilbert's Tenth Problem in its original form was to find an algorithm
to decide, given a polynomial equation $f(x_1,\dots,x_n)=0$ with
coefficients in the ring $\Z$ of integers, whether it has a solution
with $x_1,\dots,x_n \in \Z$.  In 1969 Matiyasevich~\cite{Mat70}, using
work by Davis, Putnam and Robinson (see \cite{DPR61}), proved that no
such algorithm exists, i.e.\ Hilbert's Tenth Problem is
undecidable.  Since then, analogues of this problem have been studied
by asking the same question for polynomial equations with coefficients
and solutions in other recursive commutative rings $R$.  We will refer
to this analogue of the original problem as {\sl Hilbert's Tenth
  Problem over $R$}.  Perhaps the most important unsolved problem in
this area is the case of $R=\Q$. One natural approach to showing that
Hilbert's Tenth Problem is undecidable for a ring $R$ of
characteristic 0 is to show that $\Z$ admits a Diophantine definition
over $R$, or more generally that there is a Diophantine model of the
ring $\Z$ over $R$.  We define these notions below.

\begin{definition}
\label{D:diophantine definition}
Let $R$ be a commutative ring.
Suppose $A \subseteq R^k$ for some $k \in \N$.
We say that $A$ has a {\em Diophantine definition over $R$}
if there exists a polynomial
\[
    f(t_1,\ldots,t_k, x_1,\ldots,x_n) \in R[t_1,\ldots,t_k,,x_1,\ldots,x_n]
\]
such that for any $(t_1,\ldots,t_k) \in R^k$,
\[
    (t_1,\ldots,t_k) \in A
        \quad\iff\quad
    \exists x_1,\ldots,x_n \in R, \; f(t_1,\ldots,t_k,x_1,...,x_n) = 0.
\]
In this case we also say
that $A$ is a {\em Diophantine subset} of $R^k$,
or that $A$ is {\em Diophantine over~$R$}.
\end{definition}

\begin{remark}
Suppose that $R$ is a domain whose quotient field
is not algebraically closed.
Then
\begin{enumerate}
\item[(a)]
Relaxing Definition~\ref{D:diophantine definition}
to allow an arbitrary finite conjunction of equations
in place of the single equation on the right hand side
does not enlarge the collection of Diophantine sets.
\item[(b)]
Finite unions and finite intersections of Diophantine sets are Diophantine.
\end{enumerate}
See \cite{Shbook} for details.
\end{remark}

\begin{definition}
\label{D:diophantine model}
A {\em Diophantine model of $\Z$ over a ring $R$} is a Diophantine
subset $A \subseteq R^k$ for some $k$
together with a bijection $\phi\colon \Z \to A$
such that the graphs of addition and multiplication
(subsets of $\Z^3$) correspond under $\phi$ to Diophantine subsets
of $A^3 \subseteq R^{3k}$.
\end{definition}

In 1992 Mazur formulated a conjecture that would imply that a
Diophantine definition of $\Z$ over $\Q$ does not exist, and which also
ruled out the existence of a Diophantine model of $\Z$ over $\Q$~\cite{CZ00}.  One
form of Mazur's conjecture was that for a variety $X$ over $\Q$, the
closure of $X(\Q)$ in the topological space $X(\R)$ should have at
most finitely many connected components.  This conjecture also implied
that no infinite set which is discrete in the archimedean topology has
a Diophantine definition over $\Q$.

Mazur also formulated a version of his conjecture applying to both archimedean and
nonarchimedean completions of arbitrary number fields \cite[p.\ 257]{Maz98}:
\begin{question}
\label{ques:1}
Let $V$ be any variety defined over a number field $K$.
Let $\calS$ be a finite set of places of $K$,
and consider $K_{\calS}= \prod_{v \in \calS} K_v$ viewed as
locally compact topological ring.
Let $V(K_{\calS})$ denote the topological space
of $K_{\calS}$-rational points.
For every point $p \in V(K_{\calS})$
define $W(p) \subset V$ to be the subvariety defined over $K$
that is the intersection of Zariski closures of the subsets $V(K) \cap U$,
where $U$ ranges through all open neighborhoods of $p$ in
$V(K_{\calS})$.
As $p$ ranges through the points of $V(K_{\calS})$,
are there only a finite number of distinct subvarieties $W(p)$?
\end{question}

Fix a number field $K$ and a place $\pp$.
If Question~\ref{ques:1} has a positive answer for $K$ and $\calS:=\{\pp\}$,
then there does not exist
an infinite, $\pp$-adically discrete, Diophantine subset of $K$. See
\cite[Proof of Prop.\ 1.5]{PS05} for the proof.

\bigskip
So one way to answer Question \ref{ques:1} (negatively) for $K$ would
be to
construct a Diophantine definition of
an infinite discrete $\pp$-adic set over a number field $K$.
Unfortunately, at the moment such a construction seems out of reach.
So instead we consider analogues in which $K$ is replaced by one of
its large integrally closed subrings $\OO_{K,S}$:

\begin{definition}
\label{def:bigrings}
For a number field $K$, let $\mathcal P_K$ denote the set of finite
primes of $K$, and let $\OO_{K}$ denote the ring of integers.  Given a
set~$S$ of prime ideals, not necessarily finite, the ring~$\OO_{K,S}$
is defined to be the subring of~$K$ defined by
$$\OO_{K,S}=\{x\in K:\ord_{\pp}x \geq 0  \mbox{ for all } \pp \notin S\}.$$
Observe that if $S= \emptyset$, then $\OO_{K,S}=\OO_K$ and if
$S=\calP_K$, then $\OO_{K,S}=K$.  If $S$ is finite, $\OO_{K,S}$ is
called a {\em ring of $S$-integers}. In the case where the complement
of $S$ is finite, the rings $\OO_{K,S}$
are semi-local. We will call all rings $\OO_{K,S}$ with infinite $S$
{\em big rings}.
\end{definition}

To measure the ``size'' of a set of primes one can use natural density defined below.
\begin{definition}
\label{def:natural}
Let $S \subseteq \calP_K$.
The {\em natural density} of $S$ is defined to be the limit
\[
    \lim_{X\rightarrow \infty}
    \frac{\#\{\pp \in S: N\pp \leq X\}}
    {\#\{\mbox{all } \pp: N\pp \leq X\}}
\]
if it exists.  If the limit above does not exist, one can talk about
{\it upper} density by substituting $\limsup$ for $\lim$, or {\it
  lower} density by substituting $\liminf$ for $\lim$.
\end{definition}

The study of Hilbert's Tenth Problem and of the archimedean version of
Mazur's conjecture over rings of $S$-integers has produced Diophantine
definitions of $\Z$ and discrete archimedean sets over large subrings
of some number fields (\cite{Sh97}, \cite{Sh00a},
\cite{Sh02}, \cite{Sh03},\cite{Sh04}, and \cite{Sh07}).
In 2003 Poonen proved that there exists a recursive set $S$ of primes of natural
density one such that  Hilbert's Tenth Problem is undecidable for $\Z[S^{-1}]$.
He also constructed an infinite discrete
Diophantine set (in the archimedean topology) in this ring. 
In \cite{PS05} Poonen and Shlapentokh prove
that, if there exists an elliptic curve $E$ over a number field $K$ with $\mbox{rank}(E(K))=1$, then
there exists a recursive set $S$ of primes of density one such that
Hilbert's Tenth Problem is undecidable for $\OO_{K,S}$. They also show
that there is an infinite Diophantine subset $A$ of $\OO_{K,S}$ such
that for all places $v$ of $K$, the set $A$ is discrete when viewed as
a subset of the completion $K_v$.

In \cite{EiEv09}, Eisentr\"ager and Everest reconsidered the original
result of Poonen from a different point of view, looking for a
``covering'' of $\Q$ by big rings that come from complementary sets of
primes.  More specifically, they proved that the rational primes can
be partitioned into two disjoint sets $S_1$, $S_2$ such that Hilbert's
Tenth Problem is undecidable over both $\OO_{K,S_1}$ and
$\OO_{K,S_2}$.  These results were improved by Perlega in
\cite{Perlega11} to show that the two sets can be of arbitrary
computable densities.

In this paper we generalize the results of \cite{PS05}, \cite{EiEv09} and \cite{Perlega11} to prove the following theorems:

\begin{theorem} \label{maintheorem2} Let $K$ be a number field, and
  assume there is an elliptic curve defined over $K$ with $K$-rank
  equal to 1. For every $t > 1$ and every collection
  $\delta_1,\ldots, \delta_t$ of nonnegative computable real numbers adding
  up to 1, the set of the nonarchimedean valuations of $K$ may be
  partitioned into $t$ mutually disjoint recursive subsets $S_1,\dots
  ,S_t$ of natural densities $\delta_1,\ldots, \delta_t,$ respectively,
  with the property that each ring $\OO_{K,S_i}$ contains a
  Diophantine subset discrete under any valuation of $K$ (archimedean
  or nonarchimedean).
\end{theorem}
\begin{theorem}\label{complementary}
  Assume there is an elliptic curve defined over $K$ with $K$-rank
  equal to 1. For every $t>1$ and every collection
  $\delta_1,\ldots, \delta_t$ of nonnegative computable real numbers adding
  up to 1, the set of the nonarchimedean valuations of $K$ may be
  partitioned into $t$ mutually disjoint recursive subsets $S_1,\dots
  ,S_t$ of natural densities $\delta_1,\ldots, \delta_t,$ respectively,
  with the property that $\Z$ admits a Diophantine model in each ring
  $O_{K,S_i}$. In particular, Hilbert's Tenth Problem is undecidable for
  each ring $O_{K,S_i}$.
\end{theorem}

Recently, Mazur and Rubin ~\cite{MR10} showed that if the Shafarevich-Tate
conjecture holds, then there always exists an elliptic curve defined over $K$
whose $K$-rank is one.

When proving Theorems~\ref{maintheorem2} and \ref{complementary}, we
will show that given any partition of the nonarchimedean primes into
sets $W_1,\dots, W_t$ of densities $\delta_1, \dots,\delta_t$, the
sets $S_i$ can be constructed by changing the $W_i$'s by sets of
density zero. So our results can be seen as answering the
following fundamental questions {\em up to sets of density zero}:
\pagebreak

\begin{questions}
$\left.\right.$
\be
\item For which number fields $K$ and which subsets $S$ of $\calP_K$ is
  Hilbert's Tenth Problem (un)decidable over $\OO_{K,S}$?
\item For which number fields $K$ and which subsets $S$ of $\calP_K$
is there a Diophantine model of $\Z$ over
  $\OO_{K,S}$?
\item For which number fields $K$ and subsets $S$ of $\calP_K$ is
  there an infinite subset of $\OO_{K,S}$ which is
  Diophantine over $\OO_{K,S}$ and discrete in every topology of the field $K$?

\ee
\end{questions}

One question which is not addressed by this paper is  for which number fields $K$ and which subsets $S$ of $\calP_K$  there is a Diophantine definition of $\Z$ (or $\OO_K$) over $\OO_{K,S}$.

\subsection{Overview of proof}\label{s:overview}

The goal is to prove Theorems \ref{maintheorem2} and
\ref{complementary} by partitioning $\calP_K$ into $t$ disjoint sets
$S_1, \dots, S_t$, so that each ring $\OO_{K,S_r}$ admits a
Diophantine model of the integers or has discrete infinite Diophantine
subsets. In Sections~\ref{sec:discrete} and \ref{S:construction} we
first show how to find $t$ not necessarily disjoint sets, whose union
is $\calP_K$ such that the corresponding big rings have desirable
properties. In Section~\ref{sec:complementarity} we show that these sets can also be chosen to be
mutually disjoint and of the required density.

To construct infinite discrete Diophantine sets we will proceed as in
\cite{PS05} and construct a Diophantine set containing only the elements
of a sequence converging (in all topologies of the number field) to a
limit not in the set.

To construct a Diophantine model of $\Z$ inside $\OO_{K,S_r}$, it is enough to construct a
model of the structure $$\calZ:=(\Z_{\geq 1}, 1,+, B),$$ where $B$ is a
unary predicate for the set $\{2^n + n^2: n \in
\Z_{\geq1}\}$ (see \cite[Lemma 3.16]{PS05}). A {\em Diophantine model of $\calZ$}
over a ring $R$ is a Diophantine subset $A \subseteq R^m$ for some $m$
together with a bijection $\phi: \Z_{\geq 1} \to A$ such that $\phi(B)$ is
Diophantine over $A$ and such that the graph of
addition (a subset of $\Z_{\geq 1}^3$) corresponds under $\phi $ to a
Diophantine subset of $A^3$.

In order to find suitable sets $S_r$ we  work with an elliptic
curve $E$ of rank one over $K$ and a point $P$ of infinite order that
is a suitable multiple of the generator for the non-torsion part. We
will construct $t$ (infinite) sequences of primes
\[
\{\ell_{1,1},\ell_{2,1},\dots\}, \dots, \{\ell_{1,t}, \ell_{2,t},
\dots\}
\]
such that for each $r \in \{1, \dots, t\}$, we have that
$E(\OO_{K,S_r}) \cap zE(K)$ for a suitable positive integer $z$, is
the union of $\{\, \pm \ell_{1,r} P, \pm \ell_{2,r}P, \dots\}$ and some
finite set. We then show that $A_r:=\{x_{\ell_{i,r}}: i \in
\Z_{\geq1}\}$ is a Diophantine
model of $\calZ$ in $\OO_{K,S_r}$ via the bijection $\phi: \Z_{\geq 1
} \rightarrow A_r$ sending $i$ to $x_{\ell_{i,r}}$. To prove
Theorem~\ref{maintheorem2} we construct $t$ different sequences of
primes and sets $S_r$ and show that $A_r$ as above is a discrete Diophantine set.

The paper is organized as follows. In Section~\ref{sec:recursive} we
review recursive presentations of primes of number fields, in
Section~\ref{s:primitive} we give some background about primitive
divisors and their properties, and then use these properties to prove
that certain terms in divisibility sequences have many prime ideal
divisors. Section~\ref{s:changes} describes the technical changes in
the assumptions and proofs in this paper relative to proofs and
assumptions in \cite{PS05}.  Section~\ref{sec:density} reviews and
extends some density results from \cite{PS05}.  In Sections~\ref{sec:discrete} and \ref{S:construction} we construct the rings and
the sets with the required properties.  Finally, Section~\ref{sec:complementarity} shows how to adjust the sets of primes
constructed in Sections~\ref{sec:discrete} and~\ref{S:construction} to
make them complementary.

\section{Computable Sets of Primes in Number Fields}
\label{sec:recursive}
\setcounter{equation}{0} In this section we briefly discuss a
presentation of primes in number fields and a way to define recursive
sets of primes.  We assume that a number field $K$ of degree $n$ over
$\Q$ is presented in terms of its integral basis over $\Q$.  (Such a
basis always exists and can be constructed given an irreducible
polynomial over $\Q$ of a field generator.  See for example section
7.3 of \cite{PoZa}.)  Elements of the field will be presented via
$n$-tuples of the coordinates with respect to the basis.  Given a
$K$-prime $\pp$, we will present this prime by a pair $(p,
\alpha_{\pp})$, where $p$ is the $\Q$-prime below $\pp$ and
$\alpha_{\pp} \in K$ is an algebraic integer such that
$\ord_{\pp}\alpha_{\pp}=1$ but $\ord_{\qq}\alpha_{\pp}=0$ for any
prime $\qq \not = \pp$ conjugate to $\pp$ over $\Q$.  Since the choice
$\alpha_{\pp}$ is not unique we can choose the first suitable
$\alpha_{\pp}$ under some ordering of the field.  Given an integral
basis for $K$, the map $p \mapsto
(\alpha_{\pp_1},\ldots,\alpha_{\pp_m})$, where
$p=\prod_{i=1}^k\pp_i^{e_i}$ is the factorization of $p$ in $K$, is
recursive. Further, given an element of $K$, one can effectively
determine the factorization of the divisor of this element, and given
a prime compute its norm.  Given a set of primes we can now say that
it is computable if the corresponding set of $(n+1)$-tuples
$(p,\alpha_{\pp})$ is computable.  It is also not hard to see that for
any set of $K$-primes $\calW$, the ring $\OO_{K,\calW}$ from
Definition \ref{def:bigrings} is computable if and only if $\calW$ is
computable.  For more details see Section~4 of \cite{CHS}.

\section{Primitive Divisors}\label{s:primitive}
\setcounter{equation}{0}
Let $E$ denote an elliptic
curve in Weierstrass form,
\begin{equation}\label{model}
E:\;y^2=x^3+a_4x+a_6,
\end{equation}
defined over $\OO_K$. For background, definitions and the properties of
elliptic curves used in this paper,
consult~\cite{SIL1} and~\cite{SIL2}. Let $K$ denote an
algebraic number field of degree $d=[K:\Q]$ over $\Q$.  Throughout the
paper, $E(K)$ denotes the group of $K$-rational points of $E$ and $\mathbf{O}$
denotes the point at infinity, the identity for the group of
$K$-rational points. Suppose $P$ denotes a $K$-rational point, $P\in
E(K)$, which is not torsion. Write $nP=(x_n,y_n)$. The assumptions on
$E$ allow the factorization
\begin{equation}\label{definegenbeta}
(x_n)=(x(nP))=\mathfrak a_n(P)/\mathfrak b_n^2(P)
\end{equation}
of the principal fractional ideal $(x(nP))$ into relatively prime integral ideals
$\mathfrak a_n$ and $\mathfrak b_n$.  Assuming $P$ is non-torsion guarantees that
all of the terms in the sequence $\mathfrak b=(\mathfrak b_n)$ are non-zero.

In the rational case, we may take $\mathfrak b_n$ to be a positive
integer. Silverman~\cite{SIL3} proved that when $P$ is a rational
point, for all sufficiently large
$n$, we have that $\mathfrak b_n$ has a {\it primitive divisor\/}, that is,  a
divisor of $\mathfrak b_n$ which is coprime to $\mathfrak b_m$ for all
positive integers $m<n$. In general, the expression {\it primitive
ideal divisor\/} of a term $\mathfrak b_n$ is used to describe an ideal
$\cI$ which divides $\mathfrak b_n$ but no $\mathfrak b_m$ with $m<n$. Cheon
and Hahn~\cite{ch} extended Silverman's result from~\cite{SIL3}
to algebraic number fields, showing that for all sufficiently large
$n$, it is the case that $\mathfrak b_n$ has a primitive ideal divisor.

Results about primitive divisors have a long and fine tradition for
certain sequences which satisfy a linear recurrence relation.
An interested reader can find more results concerning the existence of primitive divisors in   \cite{bhv}, \cite{epsw}, \cite{sch}, \cite{Sch}, \cite{Sch2} and \cite{Z}.\\


For a point $P$ on $E$ and a nonzero integer $n$, define $\calS _n(P)$ to be the set of all
prime ideals of $\OO_K$ that divide the ideal $\mathfrak b_n(P)$.

We will use the following properties of the sequence $\bb_n(P)$
and the sets $\mathcal{S}_n(P)$:

\begin{lemma}\label{divseq}
Let $P$ be a point of infinite order on an elliptic curve $E$ defined
over a number field $K$ as above.
\be
\item Let $n,m \in \Z - \{0\}$ and let $(m,n)$ be their gcd.  Then
  $\calS _m(P) \cap \calS _n(P) = \calS _{(m,n)}(P)$. In particular, if $P$ is an integral point of infinite order, and
  $(m,n)=1$, then $\calS _m (P)\cap \calS _n(P) = \emptyset$.
\item  The sequence $\mathfrak b_n$   is a divisibility sequence, meaning that $\mathfrak b_m   \mid\mathfrak b_n$ as ideals, whenever $m \mid n$.
\ee
\end{lemma}
\begin{proof}
$\left. \right.$
The proof follows from the standard local theory of elliptic   curves, see for example Chapters 4 and 7 in \cite{SIL1}: For $\pp
  \in \calP_K$, let $K_{\pp}$ be the completion of $K$ at $\pp$ and let
$$E_1(K_{\pp})=\{O\} \cup \{R\in E(K_{\pp}):\mbox{ord} _{\pp}(x(R))
\leq -2\}.$$ From~\cite[Proposition VII.2.2]{SIL1}, we have that $E_1(K)$ with
the elliptic curve addition is isomorphic to $\hat E(\calM)$
under the formal group addition, where $\hat E$ is the formal group
associated to $E$ and $\calM$ is the valuation ideal of $K_{\pp}$.
Note that while the proposition in \cite{SIL1} assumes that the Weierstrass equation
is minimal, this assumption is not used in the proof.  The first
assertion of the lemma now follows from the fact that $E_1(K_{\pp})$
is a group.

Furthermore, using again the fact
that the elliptic curve addition on $E_1(K_{\pp})$ corresponds to the
addition in the formal group \cite[Proposition VII.2.2]{SIL1},
together with~\cite[Corollary IV.4.4]{SIL1}, we
have that for all $R \in E_1(K_{\pp}) -\{O\}$
$$\mbox{ord} _{\pp}(x(nR)) \leq \mbox{ord} _{\pp}(x(R))-2\mbox{ord} _{\pp}(n)$$
and the second assertion of the lemma follows at once.
\end{proof}

To carry out our construction we need to prove that certain terms in
the sequence $\mathfrak b_n(P)$ have {\it many} primitive ideal
divisors. This is made precise in the next theorem.

\begin{theorem}\label{it}Let $p$ denote a prime and write $q=p^{t-1}$ for
some fixed $t\ge 2$. Suppose $Q$ is a $K$-rational point of infinite order and
$P=qQ$. Let $\{(\mathfrak b_m)(P)\}$ be the sequence of ideals coming from the
multiples of $P$ as in equation (\ref{definegenbeta}). For every large enough $n$, which is coprime to~$p$,
the term $\mathfrak b_n(P)$ has at least $t$ primitive ideal divisors. The same
is true for the terms of the sequence $\mathfrak b_{n}(pP)$.
\end{theorem}

\begin{proof}
  Let $n$ be an integer coprime to $p$ and assume that $n$ is large
  enough so that $\bb_k(Q)$ has a primitive divisor for all $k>n$. Let $\pp_{p^in}$ be
  a primitive prime ideal divisor of $\bb_{p^in}(Q)$, for $i=0,\ldots,
  t-1$.  Observe that for $i\not =j$ we have that $\pp_{p^in} \not=
  \pp_{p^jn}$.  We claim that $$\pp_{p^in} \in \calS _{p^{t-1}n}(Q)
  \setminus \calS _{p^{t-1}m} (Q)=\calS_n(P) \setminus \calS_m(P)$$ for any positive $m < n$.  Indeed,
  since $p^in$ divides $p^{t-1}n$ we have that $\pp_{p^in} \in \calS
  _{p^{t-1}n}(Q)$.  Suppose also $\pp_{p^in} \in \calS _{p^{t-1}m}(Q)$, where $m < n$.  By Lemma~\ref{divseq}, part (1),
  we can assume without loss of generality that $m$ divides $n$ and
  thus is prime to $p$.  We now also have that $\pp_{p^in} \in
  \calS _{p^{t-1}m}(Q) \cap \calS _{p^in}(Q)=\calS _{p^im}(Q)$ contradicting the assumption
  that $\pp_{p^in}$ is a primitive prime ideal divisor of
  $\bb_{p^in}(Q)$.  Thus $\pp_{p^in}, i=0,\ldots, t-1$ are primitive
  ideal divisors of $\bb_n(p^{t-1}Q)$.

  Similarly, $\pp_{p^in} \in \calS _{p^tn}(Q) \setminus \calS _{p^tm}(Q)$ for any
  positive $m < n$.  Indeed, as above, since $p^in$ divides $p^tn$ we
  have that $\pp_{p^in} \in \calS _{p^tn}(Q)$.  Suppose also $\pp_{p^in} \in
  \calS _{p^tm}(Q)$, where $0< m < n$.  Again, by Lemma~\ref{divseq},
  part (1),
  we can assume without loss of generality that $m$ divides $n$ and
  thus is prime to $p$.  We now also have that $\pp_{p^in} \in
  \calS _{p^tm}(Q) \cap \calS _{p^in}(Q)=\calS _{p^im}(Q)$, contradicting the assumption that
  $\pp_{p^in}$ is a primitive prime ideal divisor of $\bb_{p^in}(Q)$.
  Thus $\pp_{p^in}, i=0,\ldots, t-1$ are primitive ideal divisors of
  $\bb_n(p^tQ)$.

\end{proof}

\section{Some Technical Matters}\label{s:changes}
\setcounter{equation}{0}

Below we construct two collections of rings $\OO_{K,S_i}$: one to
produce infinite discrete Diophantine sets and the other to construct
a Diophantine model of the integers.  The rings $\OO_{K,S_i}$ are
constructed by generalizing the techniques from \cite{PS05}. For the
most part we use the same notation as in \cite{PS05}, but with
the following modifications:

In \cite{PS05}, the authors define ${\mathcal S}_\bad \subseteq
\calP_K$ to be the set of primes that ramify in $K/\Q$, the primes for
which the reduction of the chosen Weierstrass model is singular (this
includes all primes above $2$), and the primes at which the
coordinates of $P$ are not integral. In the rings in \cite{PS05} for
which undecidability is then shown the primes in ${\mathcal S}_\bad$
are always inverted. I.e.\ the rings are of the form $\OO_{K,S}$ with
${\mathcal S}_\bad\subseteq S$. We have to avoid inverting the primes
in ${\mathcal S}_\bad$ in each ring, otherwise the sets $S_i$ will not
be mutually disjoint. That means that in our paper the fractional
ideal generated by the $x$-coordinate of $nP$ is of the form $x(nP)=
\mathfrak a _n/\mathfrak d _n$ (with $\mathfrak a_n, \mathfrak d_n$
coprime integral ideals) and we do not have a separate ideal
$\mathfrak b _n$ that includes the contribution from the primes in
${\mathcal S}_\bad$ as in \cite{PS05}.

In view of the above, we need to show that the undecidability results
in \cite{PS05} can be proved without inverting the primes in
${\mathcal S}_\bad$. Below we note that (1) $P$ can be chosen to be
integral, that (2) we can avoid inverting the primes that ramify in
$K/\Q$ and (3) that we can avoid inverting the primes for which the
reduction of the Weierstrass model of $E$ is singular:

\begin{enumerate}
 \item We assume that the point $P:=zQ$ {\em has
  coordinates in $\OO_K$}. Here
  $Q$ generates $E(K)/E(K)_\tors$ and $z=2^{t-1}3^{t-1}\#E(K)_\tors$. This
  assumption is possible by Lem\-ma~\ref{isom} below. Our
  assumption implies that the point $P$ does not contribute any primes
  to ${\mathcal S}_\bad$.

\item Not inverting the primes that ramify in $K/ \Q$. The fact that
  ${\mathcal S}_\bad$ contains the primes that ramify in $K/\Q$ is
  used in \cite{PS05} to prove Lemma 3.3, which is then used to prove
  Proposition 3.5 in \cite{PS05}. Our proof below replaces Lemma 3.3
  and Proposition
  3.5 from \cite{PS05} with Lemma \ref{divseq} and Theorem~\ref{it}.

\item Not inverting the primes of bad reduction. Our definitions of
  $\calT_1$, $\calS_{n}$, $\mathfrak{p}_{n}, \calT_2$ differ from
  those in \cite{PS05}: Our set $\calT_1$ is contained in the set
  $\calT_1$ defined in \cite{PS05}, and it differs from it by at most
  finitely many primes (the primes in $\mathcal{S}_{\bad}$). Our set
  $\mathcal{S}_{n}$ contains {\em all} prime ideals dividing the
  denominator ideal of $x(nP)$, and $\pp_n=\mathfrak{p}_n^{(1)}$ denotes a {\em
    primitive} prime ideal divisor of the largest norm in
  $\mathcal{S}_n$. This also affects the definition of $\calT_2$.  See
  Notation~\ref{notation}  and the sets that are defined before
  Lemma~\ref{L:integer points} and Lemma~\ref{L:integer points2} below.

 The primes of bad reduction are relevant in Lemma 3.1 and
  Corollary 3.2 of \cite{PS05}. Since we have a different definition
  of $\mathfrak{p}^{(1)}_{\ell}$ we don't need to use these two results.  The only other place in \cite{PS05} where
  primes of bad reduction are relevant is Lemma 3.10, and we state
  below why this lemma still holds (see Lemmas~\ref{L:integer points}
  and \ref{L:integer points2} and their proofs).

\end{enumerate}

\begin{lemma}\label{isom}
 If $E$ is an elliptic curve and $P \in E(K)$,
then there exists a curve $E'$ that is isomorphic to $E$ over $K$ via
an isomorphism $\phi$ such that $P':=\phi(P)$ has coordinates in $\OO_K$.
\end{lemma}
\begin{proof}
If $E$ is given by a Weierstrass equation
$E:y^2 = x^3 + ax + b$ and $P \in E(K)$ has coordinates $(\alpha,\beta) \in K$,
we can choose an element $u \in K$ such that $u\alpha, u\beta \in
\OO_K$. We can then consider the curve $E'$ whose Weierstrass equation
is given by
\[
E':(y')^2=(x')^3 + au^{4}(x')+ u^{6}b,
\]
which is isomorphic to $E$ under $\phi:E\rightarrow E'$, $(x,y) \mapsto
(u^2x,u^3y)$.
The point $P':=\phi(P)$ on $E'$ has coordinates in $\OO_K$.
\end{proof}

Now we can fix some of our notation:

\subsection{Notation}\label{notation}
\setcounter{equation}{0}
\label{s:notation}
\begin{itemize}
\item Let $K$ be a number field.

\item Let $E$ be an elliptic curve of rank~1 over $K$, given by a
  Weierstrass equation with coefficients in the ring of integers $\OO_K$.
(In particular, we assume that $K$ is such that such an $E$ exists).
\item Let $E(K)_\tors$ be the torsion subgroup of $E(K)$.
\item For any set $S$ of $K$-primes let $\tilde E(\calO_{K,S})$ be the set of affine points with coordinates in $\calO_{K,S}$.
\item Let $z=2^{t-1}3^{t-1}\#E(K)_\tors$ with $t \geq 1$.
\item $P:=zQ$, where $Q$ generates $E(K)/E(K)_\tors$. As explained
  above, we may assume $P=(x,y)$ with $x,y \in \OO_K$.
\item Let $\calP_\Q=\{2,3,5,\dots\}$ be the set of rational primes.
\item Let $\calP_K$ be the set of all finite primes of $K$.

\item For $\pp \in \calP_K$, let
    \begin{itemize}
        \item $K_{\pp}$ be the completion of $K$ at $\pp$.
        \item $R_\pp$ be the valuation ring of $K_\pp$
        \item $\F_{\pp}$ be the residue field of $R_\pp$,
        \item $\norm\pp = \#\F_\pp$ be the absolute norm of $\pp$
    \end{itemize}

\item For $n \neq 0$ write $nP=(x_n,y_n)$ where $x_n,y_n \in K$.

\item Write the fractional ideal generated by $x_n$ as
\[
    (x_n)=\frac{{\mathfrak a}_n}{\mathfrak d_n},
\]
where $\mathfrak a _n$ and $\mathfrak d _n$ are coprime integral ideals.

\item For $n$ as above, let ${\calS}_n = \calS _n(P)=\{\pp \in {\mathcal P}_K :
  \pp | \mathfrak d_n\}$. By assumption on $P$, we have
  $\calS_1=\emptyset$.

\item For $\ell \in \calP_\Q$, define $a_\ell$ to be the smallest
  positive number such that ${\mathfrak d}_{\ell^{a_\ell}}$ has at
  least $t$ primitive divisors. (By Theorem~\ref{it}, applied with
  $p=2$ for $\ell \not = 2$ and with $p=3$ for $\ell=2$, we have that $a_{\ell}$ exists and $a_\ell=1$ for all but finitely many $\ell$.)
\item Let $\calL = \{\ell\in {\mathcal P}_\Q: a_\ell >1\}$
and $L = \prod_{\ell \in \calL} \ell^{a_\ell-1}$.
\item  For $k=1,\dots ,t$ define $\pp_n^{(k)}$ to
be the $k$-th largest  primitive prime divisor
of~$\mathfrak d _n$ (if it exists). (Order the primitive prime divisors according to
their norm, and break ties for prime ideals $\pp_1,\pp_2$ of the same norm according
to Section~\ref{sec:recursive}: compute the corresponding
$\alpha_{\pp_1}, \alpha_{\pp_2}$ and see which one comes first under some ordering of the field.)
\item For a prime $\ell$, define
\[
    \mu_\ell = \sup_{X \in \Z_{\ge 2}}
        \frac{\# \{\pp \in \calS_\ell : \norm\pp \le X\}}
        {\# \{\pp \in \calP_K: \norm\pp \le X\}}.
\]
\item Let ${\calM}_K$ be the set of all normalized absolute values of $K$.
\item Let ${\calM}_{K,\infty} \subset {\calM}_K$ be the set of all archimedean absolute values of $K$.
\end{itemize}

\section{On Densities of Some Sets of Primes}
\label{sec:density}
The main result of this section is the proposition below.
\begin{proposition}
\label{prop:density}
The natural density of the set $\calQ(E)=\{\qq_{\ell}, \ell \in
\Z_{>0}\}$, where $\qq_{\ell}$ is \emph{any} primitive divisor of
$[\ell]P$ (see Notation \ref{s:notation}), is zero.
\end{proposition}

In \cite{PS05} it was shown that the set $\{\pp_{\ell}, \ell \in
\Z_{>0}\}$, where $\pp_{\ell}$ is the \emph{largest} primitive divisor
of $P_{\ell}$, is equal to zero.  Below we modify this proof and show
that the primitive divisor
 does not have to be the largest in order for the density to be zero.  The key result we need from \cite{PS05} is stated below.

\begin{lemma}
For $n \in \Z_{>0}$, let $\omega(n)$ be the number of distinct prime
factors of $n$.
For any $t \ge 1$, the density of $\calZ(E,t)=\{\,\pp : \omega(\#E(\F_\pp))<t \,\}$ is $0$.  (See Lemma 3.12 of \cite{PS05}.)
\end{lemma}
As in \cite{PS05} and \cite{poonensubrings} we also need the following result and an observation.
\begin{theorem}[Hasse]
$\#E(\F_\pp) \leq \mathbf N\pp + 1 + 2\sqrt{\mathbf N \pp}.$
\end{theorem}
\begin{remark}
If $\pp$ is a prime at which $E$ has a good reduction and such that $\pp$ is a primitive divisor of $\ell P$, then $\ell | \#E(\F_\pp)$.  Note that since there are only finitely many primes at which $E$ has a bad reduction, we can ignore these primes when calculating the density.
\end{remark}
We now prove Proposition \ref{prop:density}.
\begin{proof}
We choose $\varepsilon >0$  and show that the upper natural density of $\calQ(E)$ is less than $\varepsilon$.  By the Prime Number Theorem, for some positive constants $C_{\Q}, C_K$ we have
\[
{\#\{p \in \calP_{\Q}: p \leq X\}} =O(X/\log X) < \frac{C_{\Q}X}{\log X}.
\]
\[
{\#\{\pp \in \calP_K: \mathbf N \pp \leq X\}} =O(X/\log X) > \frac{C_KX}{\log X}.
\]
Choose $t \in \Z_{>1}$ so that
\[
2^{4-t} < \frac{C_K\varepsilon}{4C_{\Q}}
\]
and choose $X \in \R_{>0}$ large enough so that
\[
\frac{ \#\{ \pp \in \calZ(E,t), {\mathbf N}\pp \leq X\}}{\#\{\pp \in \calP_K, {\mathbf N}\pp \leq X\}} < \varepsilon/2,
\]
and
\[
\frac{|\log C_K + \log\varepsilon -\log 4C_{\Q}|}{\log X} <1.
\]
Let $\overline {\calZ(E,t)}$ be the complement of $\calZ(E,t)$ in $\calP_K$.  Let $ \pp \in \overline{\calZ(E,t)}$ and assume  $\pp= \qq_{\ell}$ for some positive integer $\ell$.  In this case,
\[
\ell2^t < \#E(\F_\pp) < \mathbf N\pp + 1 + 2\sqrt{\mathbf N \pp} < 4\mathbf N \pp
\]
and therefore
\[
\ell < 2^{4-t}\mathbf N \pp \leq \frac{C_K\mathbf N \pp\varepsilon}{4C_{\Q}}.
\]
Thus for every $\pp \in \calQ(E) \cap  \overline {\calZ(E,t)}$ there exists a unique rational prime $\ell <\frac{C_K\mathbf N \pp\varepsilon}{4C_{\Q}}$.  Consider now the following ratio:
\[
\frac{\#\{\pp \in \calQ(E): \mathbf N \pp \leq X\}}{\#\{\pp \in \calP_K: \mathbf N \pp \leq X\}}=
\]
\[
\frac{\#\{\pp \in \calQ(E) \cap \calZ(E,t): \mathbf N \pp \leq X\}}{\#\{\pp \in \calP_K: \mathbf N \pp \leq X\}} + \frac{\#\{\pp \in \calQ(E) \cap \overline{ \calZ(E,t)}: \mathbf N \pp \leq X\}}{\#\{\pp \in \calP_K: \mathbf N \pp \leq X\}} \leq
\]
\[
\varepsilon/2 + \frac{\#\{\ell \in \calP_{\Q}: \ell \leq \frac{C_K\varepsilon X}{4C_{\Q}}\}}{\#\{\pp \in \calP_K: \mathbf N \pp \leq X\}} \leq
\]
\[
\varepsilon/2 + \frac{\frac{C_{\Q}C_K\varepsilon X }{4C_{\Q}\log(C_K\varepsilon X/4C_{\Q})} }{\frac{C_KX}{\log X}} =\varepsilon/2 + \frac{\varepsilon \log X}{4(\log C_K+\log\varepsilon + \log X -\log 4C_{\Q})} < \varepsilon.
\]
\end{proof}
Now we show that it is rare that $\calS_\ell$ has a large fraction of the small primes.

\begin{lemma}
\label{L:mulemma} For any $\epsilon>0$, the density of $\{\, \ell : \mu_\ell > \epsilon \,\}$ is $0$.
\end{lemma}

\begin{proof}
  The statement of this lemma is identical to the statement of Lemma
  3.8 of \cite{PS05} except for the fact that in our case $S_{\ell}$
  can contain primes of $S_{\bad}$.  However by Lemma \ref{divseq},
  only finitely many $\ell$ can be affected by the inclusion of
  $S_{\bad}$ primes and therefore the density result is unaffected.
\end{proof}
The next lemma is Lemma 3.6 of \cite{PS05} which we restate here without a proof.
\begin{lemma} \label{L:dirichlet-vinogradov}
Let $\vecalpha \in \R^n$,
let $I$ be an open neighborhood of $0$ in $\R^n/\Z^n$,
and let $d \in \Z_{\ge 1}$. Then the set of primes $\ell \equiv 1 \pmod{d}$
such that $(\ell-1)\vecalpha \mod 1$ is in $I$
has positive lower density.
\end{lemma}

\section{Infinite Diophantine Discrete Sets}
\label{sec:discrete}
\setcounter{equation}{0}
In this section we construct $t$ distinct sequences of primes
from which we will construct the sets $S_1, \dots, S_t$.
We start with a lemma which will enable us to show that the sequences we construct are computable.
\begin{lemma}
\label{computableseq}
Let $\pp_{\ell}^{(k)}$ and $\mu_{\ell}$ be as in Notation~\ref{notation}.
\be
\item For all  $k=1,\ldots,t$, the mapping  $\ell \mapsto \pp^{(k)}_{\ell}$ is computable.
\item The mapping $\ell \mapsto \mu_{\ell}$ is computable.
\ee
\end{lemma}
\begin{proof}
$\left. \right.$
\be
\item Given $k, \ell \in \Z_{> 0}$ we can effectively compute the
  coordinates of $x_{\ell}=x(\ell(P))$ and determine the factorization of
  $\dd_{\ell}$ as discussed in the Section~\ref{sec:recursive}.  By
  considering the prime factorization of $\dd_1,\ldots,\dd_{\ell-1}$
  we can determine which primes occurring in $\dd_{\ell}$ are in fact
  primitive divisors, compute their norms and determine
  $\pp_{\ell}^{(k)}$.
\item First of all, as above, for any $\ell >0$ we can effectively
  determine all the primes in $S_{\ell}$ and compute their norm.
  Secondly, once $X$ in the definition of $\mu_{\ell}$ is greater than
  the norm of $\pp_{\ell}^{(1)}$, the value of the ratio can only
  decline.  Thus to compute $\mu_{\ell}$ it sufficient to calculate
  the ratio for finitely many values of $X$ only. Therefore,
  $\mu_{\ell}$ can be computed effectively.  \ee
\end{proof}

By \cite[Corollary~VI.5.1.1]{SIL1} and \cite[Corollary~V.2.3.1]{SIL2} there is an isomorphism of real Lie groups
$\prod_{v \in \calM_{K,\infty}} E(K_v) \isom (\R/\Z)^N \times (\Z/2\Z)^{N'}$
for some $N \ge 1$ and $N' \ge 0$.
Fix such an isomorphism, and embed $E(K)$
diagonally in $\prod_{v \in \calM_{K,\infty}} E(K_v)$.
Since $P=zQ$ with $z$ even,
the point $P$ maps to an element $\vecalpha \in (\R/\Z)^{N}$.

Now we construct the sequences $\{\ell_{1,r}, \ell_{2,r}, \dots\}$ for
$r =1, \dots, t$.
To do this we describe how to define $\ell_{i,r}$ using a set $V_{i,r}, i \in
\Z_{>0}, r=1,\ldots, t$ of previously defined elements of the
sequences.  More specifically we let $V_{1,1} = \emptyset$.  For $i >
1$ we set $$V_{i,1}=\{\ell_{1,1}, \dots, \ell_{1,t},\dots,\ell_{i-1,1},
\dots, \ell_{i-1,t}\},$$ and for
$i\geq1, 1 < r \leq t$, we set
\[
V_{i,r} = \{\ell_{1,1},\dots, \ell_{1,t}, \dots,\ell_{i,1}, \dots,\ell_{i,r-1}\}.
\]
Let $\ell_{i,r}$ be the smallest prime outside $\calL$ and exceeding the bound implicit in Theorem~\ref{it} such that all of the following hold:
\begin{enumerate}
\item  $\ell_{i,r} >\ell$ for all $\ell \in V_{i,r}$,
\item \label{it:1}$\mu_{\ell_{i,r}} \le 2^{-i}$,
\item \label{it:2}$\norm\pp^{(r)}_{\ell\ell_{i,r} } > 2^i$ for all $\ell \in V_{i,r}\cup \{\ell_{i,r}\}\cup \calL$,
\item $\ell_{i,r} \equiv 1 \pmod{i!}$, and
\item $|x_{\ell_{i,r}-1}|_v > i$ for all $v \in \calM_{K,\infty}$.
\end{enumerate}
We also choose $\ell_{1,1}>3$.
\begin{proposition}
\label{P:well-defined}
The sequences $\{\ell_{1,r},\ell_{2,r},\dots\}$ are well-defined and computable for $r=1,\dots,t$.
\end{proposition}

\begin{proof}
Condition~(5) is equivalent to the requirement that
$(\ell-1)\vecalpha$ lie in a certain open neighborhood of $0$ in
$(\R/\Z)^N$, since the Lie group isomorphism maps neighborhoods of $O$
to neighborhoods of $0$.  Thus by Lemma~\ref{L:dirichlet-vinogradov},
the set of primes satisfying (4) and~(5) has positive lower
density. By Lemma~\ref{L:mulemma}, (2) fails for a set of density
$0$. Therefore it will suffice to show that (1)  and~(3) are
satisfied by all sufficiently large $\ell_i$.

For fixed $\ell$, the primes $\pp^{(r)}_{\ell\ell_{i,r} }$ for varying
values of $\ell_{i,r}$ are distinct since $\pp^{(r)}_{\ell\ell_{i,r}}$ is the $r$-th largest primitive prime divisor of $\mathfrak d_{\ell\ell_{i,r}}$.
So eventually their norms are greater than $2^i$. The same holds for $\pp^{(r)}_{\ell \ell_{i,r}}$ for fixed $\ell \in \calL$.
Thus by taking $\ell_{i,r}$ sufficiently large, we can make all the  $\pp^{(r)}_{\ell \ell_{i,r}}$ for $\ell=\ell_{i,r}$ or $\ell \in \calL$ or $\ell \in V_{i,r}$
have norm greater than $2^i$. Thus the sequence is well-defined.

Each $\ell_{i,r}$ can be computed by searching primes in increasing
order until one is found satisfying the conditions: conditions (1) --
(4) can be verified effectively by Lemma \ref{computableseq}, and
condition~(5) can be tested effectively, since $|x_{\ell_{i,r}-1}|_v$
is an algebraic real number.
\end{proof}

We now define the following subsets of $\calP_K$:
\begin{itemize}
\item $\calT_{1,r} =  \Union_{i \ge 1} \calS_{\ell_{i,r}}, r=1,\ldots,t$;
\item $\calT_{2,r}^a$ is the set of $\pp_\ell^{(r)}$
    for $\ell \notin \left(\{\ell_{1,r}, \ell_{2,r},\dots\} \cup
      \mathcal{L}\right)$, together with
    $\pp_{\ell^{a_{\ell}}}^{(r)}$ for $\ell \in \mathcal{L}$;
\item $\calT_{2,r}^b = \{\, \pp^{(r)}_{\ell_{i,r} \ell_{j,r}} : 1 \le j \le i \,\}$;
\item $\calT_{2,r}^c = \{\, \pp^{(r)}_{\ell \ell_{i,r}} : \ell \in \calL, i \ge 1 \,\}$; and
\item $\calT_{2,r} = \calT_{2,r}^a \union \calT_{2,r}^b \union \calT_{2,r}^c$.
\end{itemize}

By construction, all the terms considered above have $t$ primitive
divisors.

We now describe the important properties of these sequences.

\begin{lemma}
\label{L:integer points}
$\left.\right.$
\be
\item For each $r=1,\ldots, t$, the sets $\calT_{1,r}$ and
  $\calT_{2,r}$ are disjoint. If a subset $S_r \subset \calP_K$
  contains $\calT_{1,r}$ and is disjoint from $\calT_{2,r}$, then
  $\calE_r := \tilde E(\OO_{K, S_r}) \cap zE(K)$ is the union of $$\{\, \pm
  \ell_{i,r} P: i \ge 1 \,\}$$ and some subset of the finite set
  $\left\{\, sP : s \mid \prod_{\ell \in \calL} \ell^{a_\ell-1}
    \,\right\}$.
\item   For any $j \in \{1,2\}$ and $r,s\in \{1,\ldots,t\}$ such that
  $r \not = s$ the sets $\calT_{j,r}$ and $\calT_{j,s}$ are disjoint.
\item For any $k \in \{1,2\}$ and $r \in \{1,\ldots,t\}$ the set $\calT_{i,r}$ is computable.
\ee
\end{lemma}

\begin{proof}
$\left.\right.$
\be
\item The proof of this assertion is the same as the proof of  Lemma 3.10 of \cite{PS05}.  The proof is not affected by the fact that we do not invert primes in $S_{\bad}$, since in our case $S_1=\emptyset$ also.
\item First we assume that $j=1$.  In this case $$\calT_{1,r} \cap
  \calT_{1,s}=\left (\Union_{i \ge 1} S_{\ell_{i,r}}\right ) \cap \left
    (\Union_{i \ge 1} S_{\ell_{i,s}} \right ) =\emptyset$$ since for
  any $r\not = s$ we have that $S_{\ell_{i_1,r}} \cap S_{\ell_{i_2,s}}
  =S_{(\ell_{i_1,r},\ell_{i_2,s})}= \emptyset$ as
  $\{\ell_{i,r}\}\cap\{\ell_{i,s}\} = \emptyset$.  Next let $j=2$ and
  consider $\calT_{2,r} \cap \calT_{2,s}$.  Since the set $T_{2,r}$
  consists of the $r$-th largest primitive prime divisors of certain terms
  in the divisibility sequence $\mathfrak d_n$, and $T_{2,s}$ consists
  of the $s$-th largest primitive prime divisors of terms in the
  divisibility sequence, the definition of being a primitive divisor
  immediately implies that these sets can never have
  any nontrivial intersection when $r \neq s$.

 \item This assertion follows directly from  the fact that each
   sequence $$\{\ell_{1,r}, \ell_{2,r}, \dots\}, \, 1 \leq r \leq t$$ is computable and from Lemma~\ref{computableseq}.
\ee
\end{proof}

\begin{proposition}
\label{prop:density2}
The natural density of $\calT_{1,r}$ and $\calT_{2,r}$ $(1 \leq r \leq
t)$ is zero.
\end{proposition}
\begin{proof}
The proofs that $\calT_{1,r}$, $\calT_{2,r}^b$, and $\calT_{2,r}^c$ have density $0$
are identical to the proofs in Section~9 of~\cite{poonensubrings}.
The fact that  $\calT_{2,r}^a$ has density $0$ follows from Proposition \ref{prop:density}.
\end{proof}

Now we can construct infinite Diophantine subsets $A_r$ of
$\OO_{K,S_r}$ that are discrete in any topology of $K$.
We first need the following lemma.
\begin{lemma}\label{converge}
For each $v \in \mathcal M_K$ and $1 \leq r \leq t$ the sequence $\ell_{1,r} P, \ell_{2,r}P, \dots$
converges in $E(K_v)$ to $P$.
\end{lemma}
\begin{proof}
This is Lemma~3.14 in \cite{PS05}.
\end{proof}
We now have the following proposition.
\begin{proposition}
\label{P:Diophantine-discrete}
Let $S_r$ be as in Lemma~\ref{L:integer points} and let  $A_r:=\{x_{\ell_{1,r}},x_{\ell_{2,r} },\dots\}$.
Then $A_r$ is a Diophantine subset of $\OO_{K,S_r}$. For any $v \in \calM_K$, the set $A_r$ is discrete when viewed as a subset of $K_v$.
\end{proposition}
\begin{proof}
By Lemma~\ref{converge}, the elements of $A_r$ form a convergent
sequence in $K_v$ whose limit $x_1$ is not in $A_r$. Hence $A_r$ is
discrete.
By Lemma~\ref{L:integer points}, part (1), $x(\mathcal E_r)$ is the
union of the set $A_r$ and a finite set. Since $\mathcal E_r$ is
Diophantine over $\OO_{K,S_r}$, the set $A_r$ is Diophantine over
$\OO_{K,S_r}$ as well.
\end{proof}

In Section~\ref{sec:complementarity} we will use the sets $A_1,\ldots, A_t$ together with sets $\calT_{1,1},\ldots,\calT_{1,t}$ and $\calT_{2,1},\ldots,\calT_{2,t}$ to prove Theorem \ref{maintheorem2}.

\section{Constructing Diophantine Models of $\Z$. }
\label{S:construction}
\setcounter{equation}{0} We will now modify the $t$
sequences of primes constructed above so that in the resulting big rings Hilbert's Tenth
Problem is undecidable.

Fix two primes $\pp,\qq \in \calP_K$ of degree~1 that are primes of
good reduction for $E$, and such that $\pp$ and $\qq$ do not ramify in $K/ \Q$.
Choose $\pp,\qq$ such that neither $\pp$ nor $\qq$ divides $y_1=y(P)$,
and such that the underlying primes $p,q \in \calP_\Q$ are distinct and odd.
Let $M=p q \#E(\F_\pp) \#E(\F_\qq)$.

We now define $t$ sequences of primes $\{\ell_{i,r}\}, r=1,\ldots,t$
by using  sets
\[
V_{i,r}, i \in \Z_{>0}, r=1,\ldots, t
\]
of previously defined elements of the sequences.  More specifically we
let  $V_{1,1} = \emptyset$.  For   $i > 1$
we set $$V_{i,1}=\{\ell_{1,1}, \dots, \ell_{1,t},\dots,\ell_{i-1,1},
\dots, \ell_{i-1,t}\},$$ and for
$i\geq1, 1 < r \leq t$, we set
$$V_{i,r} = \{\ell_{1,1},\dots, \ell_{1,t}, \dots,\ell_{i,1}, \dots,\ell_{i,r-1}\}.$$
Now let $\ell_{i,r}$ be the smallest prime outside $\calL$ and exceeding the bound implicit in Theorem~\ref{it} such that all of the following hold:
\begin{enumerate}
\item $\ell_{i,r}>\ell $ for all $\ell \in V_{i,r}$,
\item $\mu_{\ell_{i,r}} \le 2^{-i}$,
\item $\norm\pp^{(r)}_{\ell\ell_{i,r}} > 2^i$ for all $\ell \in V_{i,r} \cup \{\ell_{i,r}\} \cup \calL$,
\item $\ell_{i,r} \equiv 1 \pmod M$,
\item the highest power of $p$ dividing $(\ell_{i,r}-1)/M$ is $p^i$,
and
\item $q$ divides $(\ell_{i,r}-1)/M$ if and only if $i \in B$.
\end{enumerate}

By Proposition  3.19 in \cite{PS05}  and by Proposition \ref{computableseq} we have:

\begin{proposition}
\label{P:well-defined2}
The sequences $\{\ell_{1,r},\ell_{2,r} \dots \}$, $r=1, \dots, t$, are well-defined and computable.
\end{proposition}

\begin{lemma}
\label{L:xdifference}
If $m \in \Z_{\ge 1}$, then
\[
    \ord_\pp(x_{mM+1}-x_1) = \ord_\pp(x_{M+1}-x_1) + \ord_p m.
\]
\end{lemma}

\begin{proof}
This is Lemma~3.20 in \cite{PS05}.
\end{proof}

We define the following subsets of $\calP_K$:
\begin{itemize}
\item $\calT_{1,r} =  \Union_{i \ge 1} \calS_{\ell_{i,r}}$;
\item $\calT_{2,r}^a$ is the set of $\pp_\ell^{(r)}$
    for $\ell \notin \left(\{\ell_{1,r}, \ell_{2,r},\dots\} \cup
      \mathcal{L}\right)$, together with
    $\pp_{\ell^{a_{\ell}}}^{(r)}$ for $\ell \in \mathcal{L}$;
\item $\calT_{2,r}^b = \{\, \pp^{(r)}_{\ell_{i,r} \ell_{j,r}} : 1 \le j \le i \,\}$;
\item $\calT_{2,r}^c = \{\, \pp^{(r)}_{\ell \ell_{i,r}} : \ell \in \calL, i \ge 1 \,\}$; and
\item $\calT_{2,r} = \calT_{2,r}^a \union \calT_{2,r}^b \union \calT_{2,r}^c$.
\end{itemize}

As above we now have a version of Lemma~\ref{L:integer points}.

\begin{lemma}
\label{L:integer points2}
$\left.\right.$
\be
\item For each $r=1,\ldots,t$, the sets $\calT_{1,r}$ and
  $\calT_{2,r}$ are disjoint. If a subset $S_r \subset \calP_K$
  contains $\calT_{1,r}$ and is disjoint from $\calT_{2,r}$, then
  $\calE_r := \tilde E(\calO_{K,S_r}) \cap zE(K)$ is the union of
  $\{\, \pm \ell_{i,r} P: i \ge 1 \,\}$ and some subset of the finite
  set $\left\{ \,sP : \right.$ $\left. s \mid \prod_{\ell \in \calL} \ell^{a_\ell-1}
    \,\right\}$.
\item   For any $j \in \{1,2\}$ and $r,s\in \{1,\ldots,t\}$ such that
  $r \not = s$ the sets $\calT_{j,r}$ and $\calT_{j,s}$ are disjoint.
\item For any $i \in \{1,2\}$ and $r \in \{1,\ldots,t\}$ the set $\calT_{i,r}$ is computable.
\ee
\end{lemma}

We also have an analogous version of Proposition \ref{prop:density}
and the proof is the same.

\begin{proposition}
\label{prop:density3}
The natural density of $\calT_{1,r}$ and $\calT_{2,r}$ $(1\leq r \leq t)$ is zero.
\end{proposition}

Now we can construct a Diophantine model of $\Z$ in $\OO_{K,S_r}$,
where $S_r$ is a in Lemma~\ref{L:integer points2}. We first need the
following lemma.

\begin{lemma}
  Let $B=\{\,2^n+n^2 : n \in \Z_{\ge 1}\,\}$.  Multiplication admits a
  positive existential definition in the structure $\calZ:=(\Z_{\ge
    1},1,+,B)$.  (Here $B$ is considered as a unary predicate.) Hence
  the structure $(\Z,0,1,+,\cdot)$ admits a positive existential model
  in the structure $\calZ$.
\end{lemma}
\begin{proof}
This follows from Lemma 3.16 and Corollary 3.18 in \cite{PS05}.
\end{proof}

This lemma shows that instead of finding a Diophantine model of the
ring $\Z$ over $\OO_{K,S}$, it will suffice to find a Diophantine
model of $\calZ$.

\begin{proposition}
\label{P:model of calZ}
Let $S_r$ be as in Lemma~\ref{L:integer points2} and let $A_r:=\{x_{\ell_{1,r}},x_{\ell_{2,r}},\dots\}$.
Then $A_r$ is a Diophantine model of $\calZ$ over $\OO_{K,S_r}$ via the bijection $\phi\colon \Z_{\ge 1} \to A$ taking $i$ to $x_{\ell_i,r}$.
\end{proposition}

\begin{proof}
The set $A_r$ is Diophantine over $\OO_{K,S_r}$ by part (1) of Lemma~\ref{L:integer points2}.

We have
\begin{align*}
    i \in B
    \quad&\iff\quad \text{$q$ divides $(\ell_{i,r}-1)/M$}
            &&\text{(by condition~($4$))}\\
    \quad&\iff\quad \ord_\qq(x_{\ell_{i,r}}-x_1) > \ord_\qq(x_{M+1}-x_1)
\end{align*}
by Lemma~\ref{L:xdifference} (with $\qq$ in place of $\pp$).
The latter inequality is a Diophantine condition on $x_{\ell_{i,r}}$.
Thus the subset $\phi(B)$ of $A_r$ is Diophantine over $\OO_{K,S_r}$.

Finally, for $i \in \Z_{\ge 1}$, Lemma~\ref{L:xdifference} and condition~($3$)
imply $\ord_\pp(x_{\ell_{i,r}}-x_{1}) = c + i$,
where the integer $c=\ord_\pp(x_{M+1}-x_{1})$ is independent of $i$.
Therefore, for $i,j,k \in \Z_{\ge 1}$,
we have
\[
    i+j=k
    \quad\iff\quad
    \ord_\pp(x_{\ell_{i,r}}-x_1) + \ord_\pp(x_{\ell_{j,r}}-x_1)
    = \ord_\pp(x_{\ell_{k,r}}-x_1) + c.
\]
It follows that the graph of $+$ corresponds under $\phi$
to a subset of $A_r^3$ that is Diophantine over $\OO_{K,S_r}$.

Thus $A_r$ is a Diophantine model of $\calZ$ over $\OO_{K,S_r}$.
\end{proof}

\section{Complementary Rings}
\label{sec:complementarity}
\setcounter{equation}{0} In this section we complete the proofs of
Theorems \ref{maintheorem2} and \ref{complementary}.  First we need a
general result about the existence of sets of primes of given
densities. The result we will prove is contained in
Proposition~\ref{prop:partitiondensities} below.

\subsection{Sets of primes with prescribed densities}

We start with describing the real numbers we consider as possible
densities of our sets.

\begin{definition}[Computable Real Numbers]
\label{computable reals}
A real number $\delta$ is called {\em computable} if there exists a
computable sequence of rational numbers $r_n \in \Q$ such
that $$\lim_{n\rightarrow \infty}r_n=\delta.$$
\end{definition}
It is easy to see that if $\alpha, \beta \not = 0$ are computable real numbers, then so is $\alpha/\beta$.  
Next we observe that these are the only densities we should consider in the context of our problem.
\begin{proposition}\label{prop:computable}
  Let $K$ be a number field and let $\calW_K$ be a recursive set of
  primes of $K$ having a natural density $\delta$.  In this case
  $\delta$ is a computable real number.
\end{proposition}
\begin{proof}
  Given our definition of a computable set $\calA$ of primes of a
  number field $K$ (see Section \ref{sec:recursive}), there exists a
  recursive procedure determining the size of the set $\{\pp \in
  \calA: \mathbf N \pp \leq n\}$ uniformly in $n$, and therefore there
  exists a recursive function $g(n)$ computing $r_n=
  \displaystyle\frac{\#\{\pp \in \calA: \mathbf N \pp \leq
    n\}}{\#\{\pp \in \calP_K: \mathbf N \pp \leq n\}}.$ If $\delta$ is
  the natural density of $\calA$, then $\lim_{n \rightarrow
    \infty}r_n=\delta$ and therefore $\delta$ must be a computable
  real number.

\end{proof}

We will now describe an (effective) procedure which, given a computable
positive real number $\delta$ and a computable set of primes $\calZ_K$
of a number field $K$ of natural density $\gamma \geq \delta$,
constructs a computable set of primes $\calA_K \subseteq \calZ_K$ of
natural density $\delta$.

Without loss of generality we may assume that $\gamma >
\delta$. Otherwise we set $\calA_K=\calZ_K$.  Since $\calZ_K$ is
computable, Proposition~\ref{prop:computable} implies that $\gamma$ is computable.  Let $\{d_i\},
\{g_i\}$ be computable sequences of rational numbers approximating
$\delta$ and $\gamma$ respectively. Without loss of generality
we can assume that $g_i \not =0$ for all $i \in \Z_{>0}$. Then
$\{r_i\}=\{d_i/g_i\}$ is a computable sequence of rational numbers
approximating $\alpha:= \delta/\gamma$.  Without loss of generality we may also
assume that all elements of the sequence $\{r_n\}$ are positive and strictly less
than one.

Let $\{N_i\}$ be an increasing sequence of positive
integers such that each $N_i=\mathbf N\pp_K$ for some prime $\pp_K$ of
$K$ in $\calZ_K$. Assume also that every positive integer $M \geq N_1$
occurring as a norm of a $K$-prime from $\calZ_K$ is an element of the
sequence. Let $n=[K:\Q]$ and notice that $n$ is the maximum number of
$K$-primes that can have the same norm.

We now define recursively two sequences of sets of primes
$\{\calA_i\}$ and $\{\calZ_i\}$. We denote by $a_i$ and
$z_i$ the cardinality of 
$|\calA_i|$ and $|\calZ_i|$, respectively (for $i \geq 1)$. Recall that $\{r_i\}$ is
the sequence of rational numbers (strictly less than one) approximating $\alpha=\delta/\gamma$.
\begin{enumerate}
\item Set $\calA_1=\calZ_1=\{\pp \in \calZ_K: \mathbf N\pp=N_1\}$.
\item If $\displaystyle \frac{a_i}{z_i} < r_i$, then set
  $\calA_{i+1}=\calA_i \cup \{\pp \in \calZ_K: \mathbf
  N\pp=N_{i+1}\}$.  
  Otherwise, set $\calA_{i+1} =\calA_i$.
\item For all $i \in \Z_{>0}$ set 
\[
\calZ_{i+1} = \{\pp \in \calZ_K: N_1 \leq \mathbf N\pp\leq N_{i+1}\}.
\]
\end{enumerate}
We have $\calZ_K =\bigcup_{i \in \Z_{>0}}\calZ_i$. We will now prove that the set $\calA_K := \bigcup_{i \in \Z_{>0}}\calA_i$ has
density $\delta$. To do this we need several lemmas.  It is clear from
the construction that $\calA_K$ is a computable set of primes.
\begin{lemma}
\label{le:finiteturns}$\left.\right.$
\begin{enumerate}
\item If $\displaystyle\frac{a_i}{z_i} \geq r_i$ for some index $i$,
  then there is a positive integer $k$ such that $\displaystyle\frac{a_{i+k}}{z_{i+k}} < r_i$.
\item If $\displaystyle\frac{a_i}{z_i} < r_i$ for some index $i$, then
  $\displaystyle\frac{a_{i+k}}{z_{i+k}} \geq r_i$ for some positive
  integer $k$.
\end{enumerate}
\end{lemma}
\begin{proof}
  The proofs for both assertions of the lemma will proceed by
  contradiction.  Assume that the first assertion is false, i.e.\ assume that there exists an $i \in \Z_{>0}$ such that
  $\displaystyle\frac{a_i}{z_i} \geq r_i$ and for all $k \in \Z_{>0}$
  we have $\displaystyle\frac{a_{i+k}}{z_{i+k}} \geq
  r_{i+k}$.  By step (2) in the above construction, the last inequality implies that
  $\calA_{i+k} = \calA_i$ and $a_{i+k}=a_i$ for all $k \in \Z_{>0}$.
  However, $\displaystyle \lim_{k \rightarrow \infty}z_{i+k}= \infty$
  and therefore $\displaystyle \lim_{k \rightarrow
    \infty}\frac{a_{i+k}}{z_{i+k}}=0$, while $\displaystyle \lim_{k
    \rightarrow \infty}r_{i+k} >0$, which is a contradiction. Hence
  the first assertion of the lemma is true.

Assume now that the second assertion of the lemma is false, i.e.\ that there exists an $i \in \Z_{>0}$
such that $\displaystyle\frac{a_i}{z_i} < r_i$ and such that 
\begin{equation}
\label{eq:allk}
\forall k \in \Z_{>0} \, \,\frac{a_{i+k}}{z_{i+k}} < r_{i+k}.
\end{equation}
 By step (2) of the above construction, \eqref{eq:allk} implies that 
  \[
 \forall k \in \Z_{>0} \,\, \calA_{i+k} = \calA_{i+k-1} \cup \{\pp \in \calZ_K: \mathbf N\pp=N_{i+1}\}.
  \]
This implies that $a_{i+k}\geq z_{i+k} -c$ for all $k \in \Z_{>0}$ for
some fixed nonnegative
integer $c$. At the same time, since the $r_j$'s are less than 1,
\eqref{eq:allk} implies that   $z_{i+k} > a_{i+k}$ for all $k \in \Z_{>0}$.   Thus 
\[
\lim_{k \rightarrow \infty}\frac{a_{i+k}}{z_{i+k}}=\lim_{k \rightarrow \infty}\frac{z_{i+k}}{z_{i+k}}=1 > \alpha = \lim_{k \rightarrow \infty} r_{i+k},
\]
and therefore \eqref{eq:allk} cannot hold.  
\end{proof}

We now define two sequences of positive integers that we will use below.
\begin{notation}
$\left.\right.$
\begin{itemize}
\item Let $j_1=1$ and for $i \in \Z_{>1}$ define $j_i$ to be the
  smallest positive integer greater than $j_{i-1}$ such that
  $\frac{a_{j_i\pm1}}{z_{j_i\pm1}} <\frac{a_{j_i}}{z_{j_i}}$.  (In
  other words, $\frac{a_{j_i}}{z_{j_i}}$ is a ``local
  maximum\rq{}\rq{}.)
\item Let $k_0=0$ and for $i \in \Z_{>0}$ define $k_i$ to be the
  smallest positive integer greater than $k_{i-1}$ such that 
  $\frac{a_{k_i\pm1}}{z_{k_i\pm1}} >\frac{a_{k_i}}{z_{k_i}}$.  (In
  other words, $\frac{a_{k_i}}{z_{k_i}}$ is a ``local
  minimum\rq{}\rq{}.)
\end{itemize}
\end{notation}
\begin{remark}\label{rem:props}
By construction of the sets $\calA_i,\calZ_i$ we have $z_i >a_i$ for all $i>1$. Hence if $\calA_i
  \subsetneq \calA_{i+1}$, then $\frac{a_{i+1}}{z_{i+1}}=\frac{a_i +
    m}{z_i + m} >\frac{a_i}{z_i}$, where $m$ is the number of primes
  of norm $N_{i+1}$ in $\calZ_K$.  On the other hand, if $\calA_i =
  \calA_{i+1}$ then $\frac{a_{i+1}}{z_{i+1}}=\frac{a_i}{z_i + m}
  <\frac{a_i}{z_i}$.  From Lemma \ref{le:finiteturns} we can conclude
  that both $j_i$ and $k_i$ are defined for all $i \in \Z_{>0}$.

By construction of the sets $\calA_i$ we have
  $\frac{a_{j_i-1}}{z_{j_i-1}} < r_{j_i-1}$ and
  $\frac{a_{j_i}}{z_{j_i}} \geq r_{j_i}$ for all $i >1$.  Similarly,
  we have
  $\frac{a_{k_i-1}}{z_{k_i-1}} \geq r_{k_i-1}$ and
  $\frac{a_{k_i}}{z_{k_i}} < r_{k_i}$ for all $i>0$.
\end{remark}
  We now show some properties of the sequences $\{k_i\}$ and $\{j_i\}$.

\begin{lemma}
\label{le:positions}
$\left.\right.$
\begin{enumerate}
\item For all $i \in \Z_{\geq 0}$ we have $k_i <j_{i+1} <k_{i+1}$.
\item For all $\ell \in \Z_{>0}$ there exists $i \in \Z_{>0}$ such that either 
\[
j_i \leq \ell \leq k_i \mbox{ and } \frac{a_{j_i}}{z_{j_i}} \geq \frac{a_{\ell}}{z_{\ell}} \geq \frac{a_{k_i}}{z_{k_i}}
\]
 or 
 \[
 k_i \leq \ell \leq j_{i+1} \mbox{ and } \frac{a_{k_i}}{z_{k_i}} \leq \frac{a_{\ell}}{z_{\ell}} \leq \frac{a_{j_{i+1}}}{z_{j_{i+1}}}.
 \]
\end{enumerate}
\end{lemma}
\begin{proof}
By Lemma \ref{le:finiteturns} and Remark \ref{rem:props} maxima and minima  alternate in the sequence $\{a_i/r_i\}$.  Further, by definition of $\calA_1=\calZ_1$, it is clear that $j_1=1$ produces a local maximum in the sequence.
\end{proof}
We now show that the local maxima and minima converge to $\alpha =
\delta / \gamma =
\lim_{i \rightarrow \infty}r_i $.
\begin{lemma}
\label{le:local}
$\lim_{i \rightarrow \infty}\frac{a_{j_i}}{z_{j_i}}=\alpha$ and $\lim_{i \rightarrow \infty}\frac{a_{k_i}}{z_{k_i}}=\alpha$.

\end{lemma}
\begin{proof}
  We show that for any $\varepsilon >0$ there exists a positive
  integer $M$ such that for $i >M$ we have that $|\alpha -
  \frac{a_{j_i}}{z_{j_i}}|<\varepsilon$.  The proof of the analogous
  statement with $k_i$ substituted for $j_i$ is similar. Fix $\mu < \varepsilon/4$. Let $I \in
  \Z_{>0}$ be large enough so that for all integers $s>I$ we have
  $|r_s-\alpha| <\mu < \epsilon/4$ and $z_s
  >\frac{4n(1-\alpha)}{\epsilon}$.  (Recall that
  $n=[K:\Q]$.) Fix a positive integer $s >I$ and pick an $i$ such that
  $j_i > s+1$.  By Remark~\ref{rem:props}
\begin{equation}
\label{eq:bigger}
\frac{a_{j_i}}{z_{j_i}}=\frac{a_{j_i-1}+m}{z_{j_i-1}+m} \geq r_{j_i} >\alpha -\mu,
\end{equation}
where as above, $m$ is the number of $K$-primes in $\calZ_K$ with
$K$-norm equal to $N_{j_i}$. By Remark~\ref{rem:props} we also have $\frac{a_{j_i-1}}{z_{j_i-1}} <
r_{j_i-1} < \alpha +\mu$.  Thus
\begin{align}
\label{eq:begin}\frac{a_{j_i}}{z_{j_i}} -\alpha =\frac{a_{j_i}}{z_{j_i}}-\frac{a_{j_i-1}}{z_{j_i-1}}+\frac{a_{j_i-1}}{z_{j_i-1}} -\alpha =\\
\frac{a_{j_i-1}z_{j_i-1}+mz_{j_i-1}-a_{j_i-1}z_{j_i-1}-ma_{j_i-1}}{(z_{j_i-1}+m)z_{j_i-1}}+\frac{a_{j_i-1}}{z_{j_i-1}} -\alpha =\\
\frac{m(z_{j_i-1}-a_{j_i-1})}{(z_{j_i-1}+m)z_{j_i-1}}+\frac{a_{j_i-1}}{z_{j_i-1}} -\alpha =\\
\frac{m(1-\frac{a_{j_i-1}+m}{z_{j_i-1}+m})}{z_{j_i-1}}+\frac{a_{j_i-1}}{z_{j_i-1}} -\alpha <\\
\label{eq:end}\frac{m(1-\alpha+\mu)}{z_{j_i-1}} +\mu \leq\frac{n(1-\alpha+\mu)}{z_{j_i-1}} +\mu < \varepsilon.
\end{align}
Combining \eqref{eq:bigger} and \eqref{eq:begin}--\eqref{eq:end} we conclude that
\[
\left|\frac{a_{j_i}}{z_{j_i}}-\alpha\right|<\varepsilon.
\]
\end{proof}
Now we can prove that the natural density of $\calA_K$ is $\delta$.
\begin{corollary}
The natural density of $\calA_K$ is $\delta$, i.e.\ 
$$\displaystyle\lim_{X \rightarrow \infty}\frac{\#\{\pp_K \in
    \calA_K, \mathbf N\pp_K \leq X\}}{\#\{\pp_K \in \calP_K, \mathbf
    N\pp_K \leq X\}}=\delta.$$
\end{corollary}
\begin{proof}
From Lemma \ref{le:positions} and Lemma \ref{le:local} it follows that
$\displaystyle \lim_{i \rightarrow
    \infty}\frac{a_i}{z_i}=\alpha$.  Now let $X$ be a positive real number
  greater than $N_1$.  Then $$N_i \leq X \leq N_{i+1}$$ for some
  positive integer $i$. (Recall that the $N_i$'s are norms of primes
  of $K$ that appeared in the construction of the set $\calA_K$.)
  Since $a_i$ and $z_i$ were defined to be the cardinality of
  $\calA_i$ and $\calZ_i$, respectively, we have
\[
\#\{\pp_K \in \calA_K, \mathbf N\pp_K \leq X\} = a_i,
\]
\[
\#\{\pp_K \in \calZ_K, \mathbf N\pp_K \leq X\} = z_i.
\]
This implies that $\displaystyle\lim_{X \rightarrow \infty}\frac{\#\{\pp_K
  \in \calA_K, \mathbf N\pp_K \leq X\}}{\#\{\pp_K \in \calZ_K, \mathbf
  N\pp_K \leq X\}}=\lim_{i \rightarrow \infty}\frac{a_i}{z_i}=\alpha$.

The statement of the corollary now follows from the fact that
\begin{align*}
  &\lim_{X \rightarrow \infty}\frac{\#\{\pp_K \in \calA_K, \mathbf
    N\pp_K \leq X\}}{\#\{\pp_K \in \calP_K, \mathbf N\pp_K \leq X\}}=\\
&\lim_{X \rightarrow \infty}\left(\frac{\#\{\pp_K \in \calA_K, \mathbf N\pp_K
  \leq X\}}{\#\{\pp_K \in \calP_K, \mathbf N\pp_K \leq
  X\}} \cdot \frac{\#\{\pp_K \in \calZ_K, \mathbf N\pp_K \leq X\}}{\#\{\pp_K
  \in \calZ_K, \mathbf N\pp_K \leq X\}}\right)=\\
&\lim_{X \rightarrow \infty}\frac{\#\{\pp_K \in \calA_K, \mathbf N\pp_K
  \leq X\}}{\#\{\pp_K \in \calZ_K, \mathbf N\pp_K \leq X\}} \cdot \lim_{X
  \rightarrow \infty} \frac{\#\{\pp_K \in \calZ_K, \mathbf N\pp_K \leq
  X\}}{\#\{\pp_K \in \calP_K, \mathbf N\pp_K \leq X\}}=\alpha\gamma=
\delta
\end{align*}
\end{proof}

Finally we have the following proposition.

\begin{proposition}
\label{prop:partitiondensities}
If $\delta_1, \ldots, \delta_t$ is a finite set of nonnegative
computable real numbers adding up to one, then there exist a partition
of $\calP_K$ into computable sets $W_1,\ldots,W_t$ of densities
$\delta_1, \ldots,\delta_t$, respectively.
\end{proposition}
\begin{proof}
  Without loss of generality we can assume that all the densities are
  positive.  We proceed in $t-1$ steps.  First set $\calZ_K = \calP_K$
  and construct a computable set $W_1$ of density $\delta_1$.  Observe
  that $\calP_K \setminus W_1$ is computable of density $1-\delta_1
  \geq \delta_2$.  Now set $\calZ_K =\calP_K \setminus W_1$ and
  construct $W_2, W_3$, etc.
\end{proof}
\begin{remark}
  The construction above shows in fact that there exists a partition
  of $\calP_K$ into sets of any densities adding to 1.  However, if we
  do not require that the densities are computable, the resulting sets of
  primes may be uncomputable.
\end{remark}
Now we can prove Theorems~\ref{maintheorem2} and \ref{complementary}.

\subsection{The proofs of Theorems \ref{maintheorem2} and \ref{complementary}.}

Let $\delta_1,\ldots,\delta_t$ be nonnegative computable real numbers adding
up to one.  Let $W_1,\ldots.W_t$ be a partition of primes of $K$,
where the natural density of each $W_i$ is $\delta_i$.  Such a
partition exists by Proposition \ref{prop:partitiondensities}. For the case of
Theorem \ref{maintheorem2}, let $\calT_{1,r}, \calT_{2,r},
r=1,\ldots,t$ be as defined as in Section~\ref{sec:discrete} and for
the case of Theorem~\ref{complementary}, let $\calT_{1,r},
\calT_{2,r}, r=1,\ldots,t$ be as defined as in
Section~\ref{S:construction}.  For $i = 1, \dots, t$ define
\[
S_i=(W_i \cup \calT_{1,i}\cup \calT_{2,j})\setminus (\calT_{2,i} \cup \bigcup_{r\not=i}\calT_{1,r}),
\]
 where $j \in \{1, \dots, t\}$ is such that $j\equiv i-1 \mod t$.  We claim the following:
\be
\item {\it The natural density of $S_i$ exists and is equal to
    $\delta_i$. }  This is true because by
  Propositions~\ref{prop:density} and \ref{prop:density2}, for any
  $i,j$ the natural density of $\calT_{i,j}$ is 0.
\item {\it Each $S_i$ contains all the primes of $\calT_{1,i}$ and
    omits the primes of $\calT_{2,i}$.}  To see that this assertion is
  true, observe that we explicitly add $\calT_{1,i}$ and remove
  $\calT_{2,i}$, and by Propositions \ref{L:integer points} and
  \ref{L:integer points2}, we have that $\calT_{2,i} \cap \calT_{2,j}
  =\emptyset$ for $i\not=j$.  Thus, adding $\calT_{2,j}$ does not
  introduce any primes of $\calT_{2,i}$ back.  Further from the same
  propositions removing $\bigcup_{r\not=i}\calT_{1,r}$ will not remove
  any primes of $\calT_{1,i}$.
\item {\it $S_1,\ldots,S_t$ are a partition of $\calP_K$.} First we
  show that $S_i \cap S_r =\emptyset$ for $ i\not = r$. Since $W_i$
  and $W_r$ are disjoint, the common elements can arise only from the
  primes which were added in, i.e.\ an intersection can arise from
\begin{equation}
\label{intersection}
(\calT_{1,i}\cup \calT_{2,j}) \cap (\calT_{1,r}\cup \calT_{2,l}),
\end{equation}
where $j\equiv i-1 \mod t$, and $l\equiv r-1 \mod t$ so that $l \not =
j$.  By construction, all the primes of $\calT_{1,r}$ are removed from
$S_i$ and all the primes of $\calT_{1,i}$ are removed from $S_r$.
Hence the only primes from \eqref{intersection} which can possibly be
in $S_i \cap S_r $ are in $\calT_{2,j} \cap \calT_{2,l}$.  This
intersection is empty, however, by Propositions~\ref{L:integer points}
and \ref{L:integer points2}.  Finally we show that
$\bigcup_{i=1}^t S_i =\calP_K$.  As above we start with the fact that
$\bigcup_{i=1}^t W_i =\calP_K$ and note that we only have to follow
the primes removed from $W_i$ in the process of constructing $S_i$:
\begin{equation*}
\label{union}
\calT_{2,i} \cup \bigcup_{r\not=i}\calT_{1,r}.
\end{equation*}
\ee We have shown in Part 1 of this proposition that for
$r=1,\ldots,t$, $\calT_{1,r} \subset S_r$ and
therefore the primes in the union $\bigcup_{r\not=i}\calT_{1,r}$ are
accounted for.  That leaves the primes of
\[
\calT_{2,i} \setminus \bigcup_{r\not=i}\calT_{1,r}=\calT_{2,i} \setminus \bigcup_{r=1}^t\calT_{1,r},
\]
where the equality holds because $\calT_{1,i} \cap
\calT_{2,i}=\emptyset $. When $S_i$ is constructed, this set is moved
to $S_j, j \equiv i+1 \mod t$ and observe that since $T_{2,i} \cap
T_{2,j} =\emptyset$, the primes of $\calT_{2,i} \setminus
\bigcup_{r=1}^t\calT_{1,r}$ are not removed from $S_j$.

 \section*{Acknowledgements}
Our thanks go to the London
Mathematical Society for the very enjoyable visit by the third
author to the University of East Anglia. Also, to whichever
barman invented the Marguerita.
We also thank Bjorn Poonen who suggested the use of sets of
nonnegative computable
real densities rather than rational densities.

\end{document}